\documentclass{amsart}

\usepackage[all]{xy}
\usepackage{amssymb}
\usepackage{amsmath}
\usepackage{graphicx}
\usepackage[english]{babel}

\theoremstyle{definition}
\newtheorem{definition}{Definition}
\newtheorem{remark}[definition]{Remark}
\newtheorem{example}[definition]{Example}
\newtheorem{acknowledgement}{Acknowledgement}

\theoremstyle{plain}
\newtheorem{lemma}[definition]{Lemma}
\newtheorem{proposition}[definition]{Proposition}

\newtheorem{theorem}[definition]{Theorem}

\def\t{\otimes}

\def\Prim{\mathrm{Prim\, }}
\def\Id{\mathrm{Id }}

\def\Im{\mathop{\rm Im}}

\newenvironment{theorr}[2]{\medskip\newline\noindent\textbf{Theorem \ref{#1}}\small{ (p.\pageref{#1})}\hspace{0.3cm}\normalsize\emph{#2}} {\medskip}
\newenvironment{theor}[1]{\medskip\newline\noindent\textbf{Theorem.}\hspace{0.2cm}\emph{#1}} {\medskip}

\newcommand{\F}[1]{F_#1\ }

\def\Prim{\operatorname{Prim\ }}

\def\Id{\mathrm{Id}}
\def\Im{\mathrm{Im\ }}

\def\arbreA {\vcenter{\xymatrix@R=3pt@C=3pt{
&& \\
&*{}\ar@{-}[ur] \ar@{-}[ul] \ar@{-}[d]     &\\
&&}}}

\def\arbreAt {\vcenter{\xymatrix@R=3pt@C=3pt{
&& \\
&& \\
&*{}\ar@{-}[uur] \ar@{-}[uu]\ar@{-}[uul] \ar@{-}[d]     &\\
&&}}}

\def\arbreAttt{\vcenter{\xymatrix@R=3pt@C=3pt{
&&&& \\
&&&&\\
&&*{}\ar@{-}[uurr] \ar@{-}[uur]\ar@{-}[uul] \ar@{-}[uull]\ar@{-}[dd]     &&\\
&&&&\\
&&&&}}}

\def\mcorol{\vcenter{\xymatrix@R=3pt@C=3pt{
&&m&& \\
&&\cdots&&\\
&&*{}\ar@{-}[uurr] \ar@{-}[uur]\ar@{-}[uul] \ar@{-}[uull]\ar@{-}[d]     &&\\
&&&&}}}

\def\arbreBA{\vcenter{\xymatrix@R=2pt@C=2pt{
&&&&\\
&&&*{}\ar@{-}[ul] & \\
&&*{}\ar@{-}[uurr] \ar@{-}[uull] \ar@{-}[d]     &&\\
&&&&}}}

\def\arbreBAt{\vcenter{\xymatrix@R=2pt@C=2pt{
&&&&\\
&&&*{}\ar@{-}[ul] & \\
&&*{}\ar@{-}[uurr] \ar@{-}[uu]\ar@{-}[uull] \ar@{-}[dd]     &&\\
&&&&\\
&&&&}}}

\def\arbreBAtt{\vcenter{\xymatrix@R=2pt@C=2pt{
&&&&\\
&&&*{}\ar@{-}[ul]\ar@{-}[u] & \\
&&*{}\ar@{-}[uurr] \ar@{-}[uu]\ar@{-}[uull] \ar@{-}[d]     &&\\
&&&&}}}

\def\arbreBAttt{\vcenter{\xymatrix@R=2pt@C=2pt{
&&&&\\
&&&*{}\ar@{-}[ul]\ar@{-}[u] & \\
&&*{}\ar@{-}[uurr] \ar@{-}[uull] \ar@{-}[dd]     &&\\
&&&&\\
&&&&}}}

\def\arbreAB{\vcenter{\xymatrix@R=2pt@C=2pt{
&&&&\\
&*{}\ar@{-}[ur] &&& \\
&&*{}\ar@{-}[uurr] \ar@{-}[uull] \ar@{-}[d]     &&\\
&&&&}}}

\def\arbreABttt{\vcenter{\xymatrix@R=2pt@C=2pt{
&&&&\\
&*{}\ar@{-}[ur]\ar@{-}[u] &&& \\
&&*{}\ar@{-}[uurr] \ar@{-}[uull] \ar@{-}[dd]     &&\\
&&&&\\
&&&&}}}

\def\arbreABt{\vcenter{\xymatrix@R=2pt@C=2pt{
&&&&\\
&*{}\ar@{-}[ur] &&& \\
&&*{}\ar@{-}[uurr] \ar@{-}[uu]\ar@{-}[uull] \ar@{-}[dd]     &&\\
&&&&\\
&&&&}}}

\def\arbreABtt{\vcenter{\xymatrix@R=2pt@C=2pt{
&&&&\\
&*{}\ar@{-}[ur] \ar@{-}[u]&&& \\
&&*{}\ar@{-}[uurr] \ar@{-}[uu]\ar@{-}[uull] \ar@{-}[dd]     &&\\
&&&&\\
&&&&}}}

\def\arbreAAt{\vcenter{\xymatrix@R=2pt@C=2pt{
&&&&\\
&&*{}\ar@{-}[ur]\ar@{-}[ul] && \\
&&*{}\ar@{-}[uurr] \ar@{-}[u]\ar@{-}[uull] \ar@{-}[dd]     &&\\
&&&&\\
&&&&}}}

\def\arbreAAtt{\vcenter{\xymatrix@R=2pt@C=2pt{
&&&&\\
&&*{}\ar@{-}[ur] \ar@{-}[u]&& \\
&&*{}\ar@{-}[uurr] \ar@{-}[uu]\ar@{-}[uull] \ar@{-}[dd]     &&\\
&&&&\\
&&&&}}}

\def\arbreBB{\vcenter{\xymatrix@R=2pt@C=2pt{
&&&&\\
&&&& \\
&&*{}\ar@{-}[uurr] \ar@{-}[uull] \ar@{-}[d] \ar@{-}[uu]     &&\\
&&&&}}}

\def\arbreABC{\vcenter{\xymatrix@R=1pt@C=1pt{
&&&&&&\\
&*{}\ar@{-}[ur] &&&&& \\
&&*{}\ar@{-}[uurr] &&&&\\
&&&*{}\ar@{-}[uuurrr] \ar@{-}[uuulll] \ar@{-}[d] &&&\\
&&&&&&}}}

\def\arbreBAC{\vcenter{\xymatrix@R=1pt@C=1pt{
&&&&&&\\
&&&*{}\ar@{-}[ul] &&& \\
&&*{}\ar@{-}[uurr] &&&&\\
&&&*{}\ar@{-}[uuurrr] \ar@{-}[uuulll] \ar@{-}[d] &&&\\
&&&&&&}}}

\def\arbreACA{\vcenter{\xymatrix@R=1pt@C=1pt{
&&&&&&\\
&*{}\ar@{-}[ur] &&&&*{}\ar@{-}[ul] & \\
&&&&&&\\
&&&*{}\ar@{-}[uuurrr] \ar@{-}[uuulll] \ar@{-}[d] &&&\\
&&&&&&}}}

\def\arbreCAB{\vcenter{\xymatrix@R=1pt@C=1pt{
&&&&&&\\
&&&*{}\ar@{-}[ur] &&& \\
&&&&*{}\ar@{-}[uull] &&\\
&&&*{}\ar@{-}[uuurrr] \ar@{-}[uuulll] \ar@{-}[d] &&&\\
&&&&&&}}}

\def\arbreCBA{\vcenter{\xymatrix@R=1pt@C=1pt{
&&&&&&\\
&&&&&*{}\ar@{-}[ul] & \\
&&&&*{}\ar@{-}[uull] &&\\
&&&*{}\ar@{-}[uuurrr] \ar@{-}[uuulll] \ar@{-}[d] &&&\\
&&&&&&}}}

\def\arbreCBAt{\vcenter{\xymatrix@R=1pt@C=1pt{
&&&&&&\\
&&&&&*{}\ar@{-}[ul]\ar@{-}[u] & \\
&&&&*{}\ar@{-}[uull]\ar@{-}[uu] &&\\
&&&*{}\ar@{-}[uuurrr] \ar@{-}[uuulll] \ar@{-}[d] &&&\\
&&&&&&}}}

\def\stage{\vcenter{\xymatrix@R=1pt@C=1pt{
&&&&\\
&*{}\ar@{-}[ul]&\cdots &*{}\ar@{-}[ur]&n\\
*{}\ar@{--}[rrrr]&&&&\\
&&\cdots&&\\
*{}\ar@{--}[rrrr]&&&&\\
&&&& \\
&&*{}\ar@{-}[ur] \ar@{-}[ul] \ar@{-}[d]     &&\\
&&&&1 \ .
}}}

\def\concattt{\vcenter{\xymatrix@R=1pt@C=1pt{
t_1&&t_2\\
*{}\ar@{-}[u]&&*{}\ar@{-}[u]\\
&*{}\ar@{-}[ul] \ar@{-}[dd] \ar@{-}[ur]&\\
&&\\
&&}}}

\def\concat{\vcenter{\xymatrix@R=1pt@C=1pt{
t&&s\\
*{}\ar@{-}[u]&&*{}\ar@{-}[u]\\
&*{}\ar@{-}[ul] \ar@{-}[dd] \ar@{-}[ur]&\\
&&\\
&&}}}

\def\concatt{\vcenter{\xymatrix@R=1pt@C=1pt{
t&s&u\\
*{}\ar@{-}[u]&&*{}\ar@{-}[u]\\
&*{}\ar@{-}[ul] \ar@{-}[uu]\ar@{-}[dd] \ar@{-}[ur]&\\
&&\\
&&}}}

\def\deconcattwo{\vcenter{\xymatrix@R=1pt@C=1pt{
t_1&&t_2\\
&\otimes&\\
&&\\
*{}\ar@{-}[uuu]&&*{}\ar@{-}[uuu]\\
}}}

\def\deconcat{\vcenter{\xymatrix@R=1pt@C=1pt{
t_1&&&&t_n\\&\otimes&\cdots&\otimes&\\
&&&&\\
*{}\ar@{-}[uuu]&&&&*{}\ar@{-}[uuu]\\
}}}

\def\arbreetiquette{\xymatrix@R=1pt@C=1pt{
v_{1}\ar@{-}[dd] &v_{2}\ar@{-}[dd] &...&&v_{n}\ar@{-}[dd] \\
&&&&&\\
&&&&&\\
&&...&&\\
&&&&&\\
&&&&&\\
&&*{}\ar@{-}[uurr] \ar@{-}[uull] \ar@{-}[d] &&&\\
&&&&& \ .
}}

\def\colibre{\xymatrix@R=1pt@C=1pt{
&&&&&&&\\
C \ar@{>}[ddddrrrrrrr]_{\phi} \ar@{-->}[rrrrrr]^{\tilde \phi}&&&&&&&Mag(V)^\infty=\oplus_{n \geq 0} Mag_{n}\otimes V^{\otimes n} \ar@{>>}[dddd]\\
&&&&&&&\\
&&&&&&&\\
&&&&&&&\\
&&&&&&&V \ .
}}

\def\concatcolor{\xymatrix@R=1pt@C=1pt{
T&&W\\
*{}\ar@{-}[u]&&*{}\ar@{-}[u]\\
&x{}\ar@{-}[ul] \ar@{-}[dd] \ar@{-}[ur]& \\
&&\\
&&}}

\def\produit{\xymatrix@R=1pt@C=1pt{
A^{\otimes n}\ar@{->}[dddd]_{\mu_n}&&&&&&&&\ar@{->}[llllllll]_-{\Id\otimes\cdots\otimes u\otimes\cdots\otimes\Id}A^{\otimes i}\otimes \mathbb K \otimes A^{\otimes n-i-1}\ar@{=}[rr]&&\ar@{->}[ddddllllllllll]^-{\mu_{n-1}}A^{\otimes n-1}\\
&&&&&&&&\\
&&&&&&&&\\
&&&&&&&&\\
A &&&&&&&&\\
&&&&&&&&&&&\ .
}}

\def\coproduit{\xymatrix@R=1pt@C=1pt{
C^{\otimes n}\ar@{->}[rrrrrrrr]^-{\Id\otimes\cdots\otimes c\otimes\cdots\otimes\Id}&&&&&&&&C^{\otimes i}\otimes \mathbb K \otimes C^{\otimes n-i-1}\ar@{=}[rr]&&C^{\otimes n-1}\\
&&&&&&&&\\
&&&&&&&&\\
&&&&&&&\\
C\ar@{->}[uuuu]_{\Delta_n}\ar@{->}[uuuurrrrrrrrrr]_{\Delta_{n-1}}&&&&&&&\\
&&&&&&&&&&& \ .
}}

\begin{document}

\author[E.Burgunder]{Emily Burgunder}
\address{
Institut de Math\'ematiques et de mod\'elisation de Montpellier \\
UMR CNRS 5149\\
D\'epartement de math\'ematiques\\
Universit\'e Montpellier II\\
Place Eug\`ene Bataillon\\
34095 Montpellier CEDEX 5\\
France}
\email{burgunder@math.univ-montp2.fr}
\urladdr{www.math.univ-montp2.fr/{$\sim$}burgunder/}

\title{Infinite magmatic
 bialgebras}
\subjclass[2000]{}
\keywords{ Bialgebra, Hopf algebra, Cartier-Milnor-Moore, 
Poincar\'e-Birkhoff-Witt, operad}


\begin{abstract} 
An infinite magmatic bialgebra is a vector space endowed with an $n$-ary operation, and an $n$-ary cooperation, for each $n$, verifying some compatibility relations. We prove a rigidity theorem, analogue to the Hopf-Borel theorem for commutative bialgebras: any connected 
infinite magmatic bialgebra is of the form $Mag^\infty(\Prim \mathcal{H})$, where $Mag^\infty(V)$ is the free infinite magmatic algebra over
the vector space  V.
\end{abstract}

\maketitle

\section{Introduction}
The Hopf-Borel theorem is a rigidity theorem for connected bialgebras which are both commutative and cocommutative. It takes the following form in the non-graded case:
\begin{theor}{(Hopf-Borel)
Let $\mathcal H$ be a commutative and cocommutative bialgebra,
over a field $\mathbb K$ of characteristic zero. The following are equivalent:
\begin{enumerate}
\item $\mathcal H$ is connected,
\item $\mathcal H$ is isomorphic to $S(\Prim \mathcal H)$.
\end{enumerate}}
Here $S(V)$ is the symmetric algebra over the vector space $V$, which can also be seen as the polynomial algebra.
\end{theor}

This theorem has already been generalised to other types of bialgebras, see for example \cite{L}, \cite{R} . A particular type of bialgebras, verifying a theorem analogue to the Hopf-Borel one, are magmatic bialgebras, see \cite{B}. They are vector spaces endowed with an unitary
binary operation and a counitary binary co-operation related by a magmatic compatibility relation. 
We generalise them to  bialgebras endowed with
 unitary $n$-ary operations for each $n\geq 2$, co-unitary $n$-ary co-operations , $\Delta_n$ for each $n\geq 2$, related by some infinite
magmatic compatibility relation. We denote $Mag^\infty(V)$ the free infinite magmatic algebra over a vector space $V$.

We define the primitive part of a bialgebra $H$ to be :
$$\Prim \mathcal{H}:= \cap_{n\geq 2} \left\{x \in \mathcal{H}\ \vert \  \bar \Delta_n(x)=0 \right\}\ ,$$
where,
$$\overline\Delta_n(x):=\Delta_n(x)-\sum_{i+j=n-1} 1^{\otimes i}\otimes x\otimes 1^{\otimes j}
-\sum_{m=2}^{n-1}\sum_{\sigma\in Sh_{(m,m-n)}}\sigma\circ (\overline\Delta_m(x),1^{\otimes m-n}),$$
and, $Sh_{(m,m-n)}$ are the $(m,m-n)$-shuffles.

The rigidity theorem for infinite magmatic bialgebras is as follows:
\begin{theorr}{magmatique}
{Let $\mathcal H$ be an infinite magmatic bialgebra over a field $\mathbb{K}$ of any chararacteristic.
The following are equivalent:
\begin{enumerate}
\item $\mathcal H$ is connected,
\item $\mathcal H$ is isomorphic to $Mag^\infty(\Prim \mathcal H)$.
\end{enumerate}}
\end{theorr}

The proof is based on the construction of an idempotent projector from the bialgebra to its primitive part,
 as in \cite{L}, \cite{R}, \cite{B}.

\begin{acknowledgement}
I am debtful to L. Gerritzen who raised out this question in Bochum's seminar.
I would like to thank J.-L. Loday for his advisory, 
D. Guin, R.Holtkamp for a careful reading of a first version and A.Bruguieres for some remarks. 
\end{acknowledgement}
 
\section{Infinite magmatic algebra}
\begin{definition}
An \emph{infinite magmatic algebra} $A$ is a   vector space endowed with one $n$-ary unitary operation $\mu_n$ for all $n\geq 2$ (one for each $n$)
such that:

every $\mu_n$ admits the same unit, denoted by $1$, and that,
$$
\mu_n(x_1,\cdots,x_n)=\mu_{n-1}(x_1,\cdots,x_{i-1},x_{i+1},\cdots,x_n) \quad\textrm{ where } x_i=1 \textrm{ and }  x_j\in A, \ \forall j \ . 
$$
Diagrammatically this condition is the commutativity of:
$$
\produit
$$
where $u:\mathbb K\longrightarrow A $ is the unit map.
\end{definition}
\subsection{Free infinite magmatic algebra}
\begin{definition}
An infinite magmatic algebra $A_0$ is said to be  \emph{free over the vector space $V$},
if it satisfies the following universal property. Any linear
map  $f : V\to A$, where $A$ is any infinite magmatic algebra,
extends in a unique morphism of algebras
$\tilde{f} : A_{0} \to A$:
$$
\xymatrix{
V \ar[r]^{i} \ar[dr]_{f}&A_{0} \ar[d]^{\tilde f}\\
& A& \ .}
$$
\end{definition}
\subsubsection{Planar trees and $n$-ary products}\label{pro}
A \emph{planar tree}  $T$ is a planar graph which is assumed to be simple (no loops nor multiple edges), connected and rooted. 
We denote by $Y_n$ the set of planar trees with $n$ leaves.
In low dimensions one gets:
\begin{eqnarray*}
&&Y_0=\left\{ \emptyset \right\}  ,
\ Y_1=\left\{ | \right\} ,
\ Y_2=\left\{ \arbreA \right\},
Y_3= \left\{ \arbreAB \ \arbreBA \arbreAt \right\},\\
&&Y_4=\left\{ \arbreABC \ \arbreBAC \ \arbreACA \ \arbreCAB \arbreCBA \right. \\
&&  \left.   \ \arbreABt \ \arbreBAt \ \arbreAAt \ \arbreABttt \ \arbreBAttt \ \arbreAttt \right\} , \cdots
\end{eqnarray*}
The \emph{$n$-grafting} of $n$ trees is the gluing of the  root of each tree on a new root.
For example the 2-grafting of the two trees $t$ and $s$ is:
$$
\vee_2(t,s):=
\concat \ ,
$$
the 3-grafting of three trees $t$, $s$ and $u$ is:
$$
\vee_3(t,s,u):=
\concatt \ .
$$

\begin{remark}From our definition of a planar tree, any $t\in Y_n$ is of the form
$$ t=\vee_k (t_1,\cdots,t_k)$$
for uniquely determined trees $t_1,\cdots,t_k$.
\end{remark}  

Let $V$ be a vector space. A \emph{labelled tree of degree $n$}, $n\geq 1$, denoted by $(t,v_1\cdots v_n)$,
is a tree $t$ endowed with the labelling of the leaves by the elements $v_1,\ldots, v_n$, represented as:
$$\arbreetiquette $$ 
Moreover one can define the $n$-grafting of labelled trees by the $n$-grafting of the trees, where one keeps the labellings on the leaves.

\subsubsection{Construction of the free infinite magmatic algebra}

We  denote by $Mag^\infty(V)$ the vector space spanned by the labelled planar trees:
$$Mag^\infty (V) := \oplus_{n=0}^\infty Mag_n^\infty\otimes V^{\otimes n}\ ,$$ where $Mag_n^\infty=\mathbb K[Y_n]$.

\medskip
The following result is well-known:
\begin{proposition}\label{free}
Let $V$ be a vector space.
The space $Mag^\infty(V)$ endowed with  the $n$-grafting of labelled trees, for all $n\geq 2$, is a infinite magmatic algebra. Moreover it is the free infinite magmatic algebra over $V$. 
\end{proposition}$\Box$

\section{Infinite magmatic coalgebra}
\begin{definition}
An \emph{infinite magmatic coalgebra} $C$ is a vector space endowed
 with one $n$-ary co-unitary co-operation $\Delta_n:C\to C^{\otimes n}$ for all $n\geq 2$ 
such that:

every $\Delta_n$ admits the same co-unit $c:C\longrightarrow \mathbb K$ and that the following diagram is commutative:
$$
\coproduit
$$
\end{definition}

\subsection{Construction of the connected cofree infinite magmatic coalgebra}

We denote $Sh(p,q)$ the set of $(p,q)$-shuffles. It is a permutation of $(1,\cdots,p;p+1,\cdots q)$ such that the image of the elements $1$ to $p$ and of the elements $p+1$ to $p+q$ are in order. 

We define $$\overline\Delta_n(x):=\Delta_n(x)-\sum_{i+j=n-1} 1^{\otimes i}\otimes x\otimes 1^{\otimes j}
-\sum_m\sum_{\sigma\in Sh(m,m-n)}\sigma\circ (\overline\Delta_m(x),1^{\otimes m-n}).$$
Let $T_n$ denote the $n$-corolla. Then $\overline\Delta_m(T_n)=0$ for all $m\neq n$ and $\overline\Delta_n(T_n)=|^{\otimes n}$.

\begin{definition}
An infinite magmatic co-augmented coalgebra is \emph{connected} if it verifies the following property: 
\begin{equation*}
\begin{array}{l}
\mathcal{H}= \bigcup_{r \geq 0}\F {r}\mathcal{H}\qquad
\textrm{where } \F {0}\mathcal{H}:=\mathbb{K} 1\\
\textrm{and, by induction } \F {r} \mathcal{H}:=\cap_{n\geq 2} \left\{x \in \mathcal{H}\ \vert \  \bar{ \Delta_n}(x) \in F_{r-1}\ \mathcal{H}^{\otimes n}\right\} \ ,
\end{array}
\end{equation*} 
Remark that connectedness only depends on the unit and co-operations.
\end{definition}

We define the primitive part of $\mathcal{H}$ as $\Prim \mathcal{H}:= \cap_{n\geq 2} \left\{x \in \mathcal{H}\ \vert \  \bar \Delta_n(x)=0 \right\}$.

\begin{definition}
An infinite magmatic coalgebra $C_{0}$ is \emph{cofree on the vector space $V$} if there exists a linear map 
$p :C_{0} \to V$ satisfying the following universal property:

any linear map $\phi : C\to V$, where $C$ is any connected infinite magmatic coalgebra 
such that $\phi(1)=0$,
extends in a unique coalgebra morphism $\tilde \phi : C\to C_{0}$:
$$
\xymatrix{C \ar[dr]^{\phi} \ar[d]_{\tilde \phi}&\\
C_{0} \ar[r]^{p} & V \ .}
$$
\end{definition}
\subsubsection{Planar trees and $n$-ary coproducts}\label{copro}

We endow the vector space of  planar trees with the following $n$-ary co-operations, for $n\geq 2$:
for any planar tree $t$ we define:
$$\Delta_n(t):=\sum t_1\otimes\cdots\otimes t_n$$
where the sum is extended on all the ways to write $t$ as $\vee_n(t_1,\cdots,t_n)$, where $t_i$ may be $\emptyset$.
It can be explicited, as follows,  
for $t=\vee_n(t_1,\cdots,t_n)$, where $t_i\neq \emptyset$ for all $i$: 
\begin{eqnarray*}
&&\Delta_n(t):=
\left(\deconcat \right)+\sum_{i=0}^{n-1}\emptyset^{\otimes i}\otimes t\otimes \emptyset^{\otimes n-i-1}\ ,\\
&&\Delta_m(t):=\left\{ \begin{array}{ll}
\sum_{i=0}^{m-1}\emptyset^{\otimes i}\otimes t\otimes \emptyset^{\otimes m-i-1}\ ,
& \textrm{if $m<n$}\\
\\
\sum_{i=0}^{m-1}\emptyset^{\otimes i}\otimes t\otimes \emptyset^{\otimes m-i-1}+\\
+\sum_{i_1+\cdots+i_{n+1}=m-n}\emptyset^{\otimes i_1}\otimes t_1 \otimes \emptyset^{\otimes i_2}\otimes \cdots\otimes t_n\otimes \emptyset^{\otimes i_{n+1}}
\ ,
& \textrm{if $m>n$}\end{array} \right.\\
&&\Delta_n(|):=\sum_{i=0}^{n-1} \emptyset^{\otimes i}\otimes |\otimes \emptyset^{n-i-1} \ ,\\
&&\Delta_n(\emptyset):=\emptyset^{\otimes n} \ .
\end{eqnarray*}
As in the preceding section one can define the $n$-ungrafting of labelled trees by
 the $n$-ungrafting of planar trees and keeping the labelling on the leaves.
 
 Remark that the empty tree $\emptyset$ plays here the role of the unit, it can then be denoted by $1:=\emptyset$.
\subsubsection{Construction of the cofree connected infinite magmatic coalgebra}

\begin{definition}
The \emph{height} of a planar tree $T$ is the maximal number of inner vertices one can meet when going through all the paths starting from 
the root to a leaf.
$$\stage$$ 
\end{definition}

\begin{example}
The $n$-corolla is of height 1.
The tree $\arbreCBAt$ is of height $3$.
\end{example}

\begin{proposition}
Let $V$ be a vector space.
The space $Mag^\infty(V)$ endowed with the $n$-ungrafting co-operations on labelled trees
 is a connected infinite magmatic coalgebra. Moreover it is cofree 
over $V$  among the connected infinite magmatic coalgebras.
\end{proposition}
\begin{proof}
We could prove this proposition by dualising Proposition (\ref{free}), but since we did not 
give a proof of it we will write completely this proof.

The co-operations are co-unital by definition, so $Mag^\infty$ is an magmatic coalgebra.

Then we  verify the connectedness of $Mag^\infty(V)$.
It comes naturally that:
\begin{eqnarray*}
\F 1  Mag^\infty (V) &=& Mag_0(V)\oplus Mag_1(V)\\
\F 2  Mag^\infty (V) &=& Mag_0(V)\oplus Mag_1(V) \oplus_n \{n- \textrm{corollas}\}
\end{eqnarray*}
One can conclude by induction on the number of heights of the tree.
Indeed, let us consider the tree $T\in Mag^\infty(V)$. It can
 be seen as the $n$-grafting of other trees, each of them having at least a height
less than the considered tree. Moreover we have:
\begin{eqnarray*}
\overline\Delta_n(T)&=&\Delta_n\circ\mu_m(T_1\otimes\cdots\otimes T_n)\\ 
&=&\left\{ \begin{array}{ll}
0 \ ,
& \textrm{if $m\neq n$}\\
\\
\underbrace{T_1}_{\in \F i_1(V)}\otimes\cdots\otimes \underbrace{T_n}_{\in \F i_n(V)} \ , 
& \textrm{if $m=n$.}\end{array} \right.
\end{eqnarray*}
where $i_1,\cdots,i_n \leq n$, so $ \overline\Delta_n (T) \in \F {n-1}^{\otimes n}$.
So we can conclude that:
\begin{equation*}
\F r Mag^\infty (V)=\oplus_{m=0}^{m=r} \{ \textrm{ trees with height } m \}\ .
\end{equation*}
It is clear that $\cup_n \F n Mag^\infty (V) = Mag^\infty(V)$. 

To prove the cofreeness of the coalgebra, it is sufficient to prove the commutativity of the
following diagram:
\begin{equation}\label{colibre}
\colibre
\end{equation}

The map $\tilde{\phi}$ can be decomposed into its homogeneous components as follows:
\begin{equation}
\tilde{\phi} (c)= \tilde{\phi} (c)_{(1)}+ \tilde{\phi} (c)_{(2)} + \tilde{\phi}(c)_{(3)}+...
\end{equation}
By induction on $n$, one can determine the homogenous components of $\tilde{\phi}$.
As the map $\tilde{\phi}$ is a coalgebra morphism defined on $\bar C$, one defines $\tilde{\phi} (1)=1$ .
 
The commutativity of the diagram (\ref{colibre}) gives the following equality:

\begin{equation}\label{egalite sur V}
\tilde{\phi}(c)_{1}=(|,\phi(c)).
\end{equation}

By definition of $Mag_{2}(V)$:

\begin{equation*}
\tilde{\phi}(c)_{2}=\sum \ (\arbreA,a_{1}a_{2})
\end{equation*}
We adopt the following notation $\bar{\Delta}(c)= \Sigma c_{1}\otimes c_{2}$. And we compute:
\begin{equation*}
\begin{array}{rcl}
\tilde{\phi}_1\otimes \tilde{\phi}_1\circ \bar{\Delta}(c) & = & \sum \tilde{\phi}_{1}(c_{1})\otimes \tilde{\phi}_{1}(c_{2})\\
&=&\sum (|, \phi (c_{1})) \otimes (|,\phi (c_{2})) \ \textrm{thanks to (\ref{egalite sur V})}\\
\multicolumn{3}{l}
{\textrm{But }
\Delta \circ \tilde{\phi}(c)_{2}=\sum (|,a_{1}) \otimes (|,a_{2})}\\
&=&(|,\phi (c_{1}))\otimes (|,\phi (c_{2}))\\
\multicolumn{3}{l}
{\textrm{Therefore, }}\\
\tilde{\phi}(c)_{2}&=&\sum (\arbreA , \phi(c_{1}) \phi (c_{2}))
\end{array}
\end{equation*}
Any tree $T$ determines a co-operation that we denote by $\Delta^T$. If $T$ is the corolla
with $n$ leaves then $\Delta^T$ is $\Delta_n$. Another example is to consider the tree  $T={\arbreCBAt}$, we have 
\begin{eqnarray*}
\Delta^T=(\Id^{\otimes 3}\otimes\Delta_3) \circ(\Id\otimes\Delta_3)\circ \Delta_2 \ .
\end{eqnarray*}
Analogously for a tree $T_i$ of degree $n$:
\begin{eqnarray*}
\tilde\phi_n(c)&=&\sum(t,a_1\cdots a_n)\\
\bar\Delta^{T_i}\tilde\phi(c)&=&\sum(|,a_1^i)\otimes\cdots\otimes(|,a_n^i)\\
\textrm{Denote: }\bar\Delta^{T_i}(c)&=&\sum c_1^i\otimes\cdots\otimes c_n^i\\
\tilde\phi_1^{\otimes n}\circ\bar\Delta^{T_i}(c)&=&\sum(|,\phi(c_1^i))\otimes\cdots\otimes(|,\phi(c_n^i)),
\end{eqnarray*}
which gives us:
$$\tilde\phi^i_n(c)=\sum(T_i,\phi(c_1^i)\cdots\phi(c_n^i)) \ .$$
Going through all the trees of degree $n$, we have:
$$\tilde\phi_n(c)=\sum_{T_i \textrm{ of degree n}}(T_i,\phi(c_1^i)\cdots\phi(c_n^i))$$
(though we denote $\tilde\phi^i , \ T_i$, we don't assume that there must be an order on the trees, this notation is only used to  distinguish the trees with same degree.)

Therefore one has:  
\begin{equation*}
\begin{array}{rcl}
\tilde{\phi}(c)&=&(|,\phi(c)) + \sum  (\arbreA ,\phi (c_{1}) \phi(c_{2}))+\\ 
&&\sum (\arbreAB,\phi (c_{1}^{1}) \phi (c_{2}^{1}) \phi (c_{3}^{1}))
+ \sum (\arbreBA ,\phi (c_{1}^{2}) \phi (c_{2}^{2}) \phi (c_{3}^{2}))+\\
&&\sum \Big(\arbreAt,\phi (c_{1}^{3}) \phi c_{2}^{3} \phi (c_{3}^{3})
\Big)+...\\
\end{array}
\end{equation*}

By construction $\tilde{\phi}$ is a morphism of infinite magmatic coalgebras  which is unique, since we  have no other choice to have the commutativity of diagram (\ref{colibre})  and  the coalgebra morphism property.
\end{proof}
\section{Infinite magmatic bialgebra}
\begin{definition}
An \emph{infinite magmatic bialgebra} $(\mathcal{H},\mu_n,\Delta_n)$ is a vector space
$\mathcal{H}=\bar {\mathcal{H}}\oplus \mathbb{K} 1$ such that:\\
1) $\mathcal{H}$ admits an infinite magmatic algebra structure with $n$-ary operations denoted~$\mu_n$\!\!~,\\
2) $\mathcal{H}$ admits a infinite magmatic coalgebra structure with $n$-ary co-operations denoted~$\Delta_n$,\\
3) $\mathcal{H}$ satisfies the following compatibility relation called the ``infinite magmatic compatibilty'':
\begin{equation}\label{rel bigebres}
\begin{array}{l}
\Delta_n\circ\mu_n (x_1\otimes\cdots\otimes x_n)=
x_1\otimes\cdots\otimes x_n+\sum_{i=0}^{n-1}1^{\otimes i}\otimes \underline x \otimes 1^{\otimes n-i-1}\ ,\\
\Delta_m\circ\mu_n (x_1\otimes\cdots\otimes x_n)=\\
\qquad\qquad\qquad\left\{ \begin{array}{ll}
\sum_{i=0}^{m-1}1^{\otimes i}\otimes \underline x \otimes 1^{\otimes m-i-1}\ ,
& \textrm{if $m<n$}\\
\\
\sum_{i=0}^{m-1}1^{\otimes i}\otimes \underline x \otimes 1^{\otimes m-i-1}+\\
+\sum_{i_1+\cdots+i_n+1=m-n}1^{\otimes i_1}\otimes x_1 \otimes 1^{\otimes i_2}\otimes \cdots\otimes x_n\otimes 1^{\otimes i_{n+1}}
& \textrm{if $m>n$}\end{array} \right.
\end{array}
\end{equation} 
where $  \underline x :=\mu_n \circ x_1\otimes\cdots\otimes x_n \ and \  x_1,\cdots,x_n \in \bar{\mathcal{H}}$ .
\end{definition}

A fundamental example in our context is the following:

\begin{proposition}
Let $V$ be a vector space.
The space
$(Mag^\infty(V),\vee_n,\Delta_n)$, where the operations $\vee_n$ (resp. the co-operations $\Delta_n$) are defined in \ref{pro} and \ref{copro}, is an
 infinite magmatic connected bialgebra.
\end{proposition}

\begin{proof}
Any tree can be seen as the $n$-grafting of $n$ trees, except the empty tree and the tree reduced to the root. Therefore the $m$-ungrafting of a tree can be viewed 	as the $m$-ungrafting of the $n$-grafting of $n$ trees. This observation gives the compatibility relation.
\end{proof}

\section{The main theorem}

\begin{definition}
The \emph{completed infinite magmatic algebra}, denoted by
$Mag^\infty(\mathbb K)^{\wedge}\ ,$ is defined by
$$Mag^\infty(\mathbb K)^{\wedge}~\!\!=~\!\!\prod_{n\geq 0}Mag_n \ ,$$ where
 the first generator  $|$ is denoted by $t$. 
This definition allows us to define
 formal power series of trees in $Mag^\infty(\mathbb K)^{\wedge}$.
\end{definition}

\begin{lemma}\label{finv}
The following two formal power series, $g$ and $f$, are inverse for composition in
$Mag^\infty(\mathbb K)^{\wedge}$:
$$g(|):=|-\arbreA-\arbreAt-\arbreAttt-\cdots , \qquad f(|):= \sum T, $$ where the sum is extended to all planar trees $T$. 
\end{lemma}
Here the tree $T$ stands for the element $T(x):=T(x,\ldots,x)$, where $x=|$ the generator. The composition of 
$T_1\circ T_2$ is defined as $T_1\circ T_2(x):=T_1\circ T_2(x,\ldots,x)=T_1(T_2(x,\ldots,x),\ldots,T_2(x,\ldots , x))$.

\begin{proof}
First, we show that $g\circ f = |$, that is to say: $$\sum T-\sum_{T_1,T_2} \vee_2 (T_1\otimes T_2)-\cdots -\sum_{T_1,\cdots, T_n}\vee_n (T_1\otimes\cdots\otimes T_n)-\cdots= | \ ,$$ equivalently:
$$\sum_{T_1,T_2} \vee_2 (T_1\otimes T_2)-\cdots -\sum_{T_1,\cdots, T_n}\vee_n (T_1\otimes\cdots\otimes T_n)-\cdots= \sum T - | \ .$$
It is immediate, as every tree can be seen as the $n$-grafting of $n$ trees for a certain $n$, except $|$.

Then one verifies that, as in the associative case, a right inverse is also a left inverse.
Let $f^{-1}$ denote the left inverse of  $f$. Then:
\begin{eqnarray*}
f^{-1}=f^{-1}\circ (f\circ g)=(f^{-1}\circ f)\circ g = g \ .
\end{eqnarray*}
Remark that we have associativity of composition even in the infinite magmatic context.
Therefore one has $f\circ g=Id$ et $g\circ f=Id$.
\end{proof}

\begin{definition}
The \emph{$n$-convolution} of $n$ infinite magmatic algebra morphisms $f_1,\cdots ,f_n$ is defined by:
\begin{equation*}
 \star_n (f_1\cdots f_n) := \mu_n \circ (f_1 \otimes \cdots \otimes f_n) \circ \Delta_n  \ .
\end{equation*}
Observe that these operations are unitary.
\end{definition}

\begin{lemma}
Let $(\mathcal{H}, \mu_n, \Delta_n)$ be a connected infinite  magmatic bialgebra. The linear map
$e:\mathcal{H} \rightarrow \mathcal{H}$ defined as: 
\begin{equation*}
e:=J- \star_2\circ J^{\otimes 2}-\star_3\circ J^{\otimes 3}-\cdots-\star_n\circ J^{\otimes n}-\cdots
\end{equation*}
where $J= Id-uc$, $u$ the unit of the operations, $c$ the co-unit of the co-operations,
has the following properties:
\begin{enumerate}
\item $\Im e = \Prim \mathcal{H} $,
\item for all $x_1,\cdots,x_n \in \bar{ \mathcal{H}}$ one has $ \ e\circ\mu_n (x_1\otimes \cdots \otimes x_n)=0$,
\item the linear map $e$ is an idempotent,
\item for $\mathcal{H}=(Mag^\infty(V),\mu_n,\Delta_n)$ defined above, $e$ is the identity on $V = Mag_1(V)$
and trivial on the other components.
\end{enumerate}
\end{lemma}

\begin{proof} 
In this proof, we adopt the following notation:
$Id:=Id_{\bar{ \mathcal{H}}},$ and for all $x \in \bar{ \mathcal{H}}, \ 
\bar \Delta_n(x):=\sum x_{1}\otimes\cdots\otimes x_{n}$ 

\begin{enumerate}
\item Proof of $\Im e = \Prim \mathcal{H} $ \ .
\begin{eqnarray*}
\overline \Delta_n (e(x))&=&\overline \Delta_n(x)-\sum_m\overline{\Delta_n}\circ\mu_m\circ\overline {\Delta_n}(x)\\
&=&x_1\otimes\cdots\otimes x_n-\overline \Delta_n\circ\mu_n(x_1\otimes\cdots\otimes x_n-\sum_{m\neq n}
\underbrace{\overline\Delta_n\circ\mu_m\circ\overline {\Delta_n}(x)}_{=0}\\
&=& 0\ .
\end{eqnarray*}
\item Proof that for all $ x_1,\cdots ,x_n \in \bar{ \mathcal{H}} $ one has $e\circ \mu_n(x_1\otimes\cdots\otimes x_n)=0$ .
Indeed,
\begin{eqnarray*}
e\circ \mu_n (x_1\otimes\cdot\otimes x_n)&=& \mu_n (x_1\otimes\cdots\otimes x_n)-\sum_m\mu_m
\circ\overline{\Delta_m}\circ\mu_n (x_1\otimes\cdots\otimes x_n)\\
&=&\mu_n (x_1\otimes\cdots\otimes x_n)-\mu_n
\circ\overline{\Delta_n}\circ\mu_n (x_1\otimes\cdots\otimes x_n)\\
&=& 0 \ .
\end{eqnarray*}

\item Proof that $e$ is an idempotent.
We compute:
\begin{eqnarray*}
e(e(x)) &=& e(x)-\sum_m e(\mu_m\circ\overline {\Delta_n}(x))\\
&=& e(x). 
\end{eqnarray*}

\item Proof that for $\mathcal{H}=(Mag^\infty(V),\mu_n,\Delta_n)$ defined above, $e$ is the identity on $V = Mag_1(V)$
and trivial on the other components.

On $Mag_{1}(V) = | \otimes V$ we have: $e(|\otimes x)= |\otimes x$.
All other trees can be seen as the $n$-grafting of $n$ trees for a certain $n$. Then it suffices to apply the second property 
of the idempotent $e$ to complete the proof.
\end{enumerate}
\end{proof}

\begin{theorem}\label{magmatique}
If  $\mathcal{H}$ be a connected infinite magmatic bialgebra over a field $\mathbb{K}$ of any characteristic, then the following are equivalent:
\begin{enumerate}
\item $\mathcal{H}$ is connected,
\item $\mathcal{H} \cong Mag^\infty(\Prim \mathcal{H})$.
\end{enumerate}
\end{theorem}

\begin{proof}
It is convenient to introduce the following notation: for $T\in Y_n$
$$\star_T (J):\mathcal H\longrightarrow \mathcal H:x\mapsto x^{n}\mapsto \star_T (J)(x^n)$$where $T\in Mag(\mathbb K)$: we label the tree $T$ by  $J$ on each leaf, and endow each inner vertex by an $n$-ary operation $\star_n$. For example considering the tree  $T={\arbreCBAt}$, we have 
$$\star_T (J) = \star_2(J)\circ(\Id\otimes \star_3(J))\circ(\Id^{\otimes 3}\otimes\star_3(J))\ , $$
and valued on an element of $\mathcal H$, one has:
$$\star_T (J)(x) = \star_2(J)\circ(\Id\otimes \star_3(J))\circ(\Id^{\otimes 3}\otimes\star_3(J))(x^{\t 6})\ , $$
Observe that $$J^{\star T_1}\otimes\cdots\otimes J^{\star T_n} = J^{\star(\vee_n( T_1 \otimes\cdots\otimes T_n))}$$ by definition.
Let us denote $V := \Prim \mathcal{H}$. 

We prove the isomorphism by explicitly giving the two inverse maps.

We define: 
$G:\bar{ \mathcal{H}} \rightarrow \overline{Mag^\infty(V)}$ as
\begin{equation*}
G(x):=J(x)-\star_2\circ J^{\otimes 2}(x)-\star_3\circ J^{\otimes 3}(x)-\cdots-\star_n\circ J^{\otimes n}(x)-\cdots \ ,
\end{equation*}
and  $F:\overline{ Mag}(V) \rightarrow \bar{ \mathcal{H}}$ by
\begin{equation*}
F(x):=\sum {\star_ T}(J)(x) \ ,
\end{equation*}
where the sum is extended to all planar trees $T$.

Moreover, denote by $t$ the generator of $Mag^\infty(\mathbb K)$, $t:=|$, and by $t^n := \vee_n\circ t^{\otimes n}$. We 
define $g(t):=t -t^{2}-t^3-\cdots-t^n-\cdots,$ and $f(t):= \sum T $, where the sum is extended to
all planar trees $T$.
By lemma ($\ref{finv}$) these two preceding maps are inverse, for composition. 

These series can be applied to elements of  $Hom_\mathbb{K}(\mathcal{H},\mathcal{H})$ sending  $1$ on $0$ using $\star_n$ as a product,
thanks to the following morphism:
\begin{equation*}
\begin{array}{rcl}
Mag(V)^{\infty\wedge}& \longrightarrow & Hom_\mathbb{K}(\mathcal{H},\mathcal{H})\\
t&\mapsto & J\\
\phi(t)=\sum a_n T&\mapsto &\phi^\star(J)=\Phi=\sum a_nJ^{\star T}\\
\phi\circ\psi(t)&\mapsto & (\phi\circ\psi)^\star(J) =\Phi\circ\Psi=\phi ^\star(J)\circ\psi ^\star(J) 
\end{array}
\end{equation*}
It is clear that $e=g^\star (J) $. 

Therefore composing the two power formal series $F$ and $G$ gives as a result: 
\begin{eqnarray*}
F \circ G =f^{\star} \circ g^{\star}(J)= (f \circ g)^{\star}(J)=Id^{\star}(J)=J\\
G \circ F =g^{\star} \circ f^{\star}(J)= (g \circ f)^{\star}(J)=Id^{\star}(J)=J
\end{eqnarray*}
The proof is complete since $J=Id$ on $\bar{ \mathcal{H}}$. \\
\end{proof}

\begin{remark}
There is a slightly different definition under which the result still holds.
Change the definition  in:

An infinite magmatic co-augmented coalgebra is \emph{connected} if it verifies the following property: 
\begin{equation*}
\begin{array}{l}
\mathcal{H}= \bigcup_{r \geq 0}\F {r}\mathcal{H}\qquad
\textrm{where } \F {0} \mathcal{H}:=\mathbb{K} 1\\
\textrm{and, by induction } \F {r} \mathcal{H}:=\cap_{n\geq 2}^{r+1} \left\{x \in \mathcal{H} \ \vert \  \bar {\Delta_n}(x) \in F_{r-1}\ \mathcal{H} ^{\otimes n}\right\} \ ,
\end{array}
\end{equation*}
The definition of the primitive elements being unchanged.
Then, we find that $Mag^\infty (\mathbb K)$ is still connected for the following description:
\begin{eqnarray*}
\F 1 Mag^\infty &=& \{ | \} \\
\F 2 Mag^\infty &=& \{ |, \arbreA \} \\
\F n Mag^\infty &=& \{ \textrm{the n corolla and all the trees with a root being a $m$-grafting} \\
&&\textrm{where $m\leq n$ and all the operations being in }
\F n Mag^\infty  \}
\end{eqnarray*} 
We observe that the space of primitive elements $\Prim \mathcal{H}$ is the same as defined in the precedent cas, that is to say $\Prim \mathcal{H}=\{|\}.$
\end{remark}

\section{$m$-magmatic bialgebras}

Instead of considering infinite magmatic bialgebras one may consider $m$-magmatic bialgebras, with $m\geq 2$, where the number of operations and co-operations is restricted to $m$.
Explicitly, we would have:

\subsection{ $m$-magmatic algebra and free $m$-magmatic algebra}

\begin{definition}
An \emph{$m$-magmatic algebra} $A$ is a   vector space endowed with one $n$-ary unitary operation $\mu_n$ for all $2\leq n\leq m$ (one for each $n$)
such that:

every $\mu_n$ admits the same unit $u$ and that,
$$
\mu_n(x_1,\cdots,x_n)=\mu_{n-1}(x_1,\cdots,x_{i-1},x_{i+1},\cdots,x_n) \quad\textrm{ where } x_i=1 \textrm{ and }  x_j\in A, \ \forall j \ . 
$$
Diagrammatically this condition is the commutativity of:
$$
\produit
$$
\end{definition}

\subsubsection{Construction of the free $m$-magmatic algebra over a vector space $V$}\label{mpr}

Let us describe the set of $m$-ary planar trees.
A \emph{$m$-ary planar tree}  $T$ is a planar graph which is assumed to be simple (no loops nor multiple edges) and connected, such that the valence of each inner vertex is at most $m+1$. 
We denote by $Y_n^m$ the set of $m$-ary planar trees with $n$ leaves.
In low dimensions one gets:
\begin{eqnarray*}
&&Y_0^m=\left\{ \emptyset \right\}  ,
\ Y_1^m=\left\{ | \right\} ,
\ Y_2^m=\left\{ \arbreA \right\},\cdots
\end{eqnarray*}
For example if $m=2$ we have the binary planar trees, see \cite{B}.
For $m=3$, in low dimensions we  have:
\begin{eqnarray*}&&Y_3^3= \left\{ \arbreAB \ \arbreBA \arbreAt \right\},\\
&&Y_4^3=\left\{ \arbreABC \ \arbreBAC \ \arbreACA \ \arbreCAB \arbreCBA \right. \\
&&  \left.   \ \arbreABt \ \arbreBAt \ \arbreAAt \ \arbreABttt \ \arbreBAttt  \right\} , \cdots
\end{eqnarray*}
We observe that we have lost the tree $\arbreAttt$ from the planar case.

The \emph{$n$-grafting}, $2\leq n\leq m$, of $n$ trees is the gluing of the  root of each tree on a new root, exactly as in the infinite magmatic case.

\begin{remark}From our definition of a $m$-ary planar tree, any $t\in Y_n^m$ is of the form
$$ t=\vee_k (t_1,\cdots,t_k)$$
for uniquely determined trees $t_1,\cdots,t_n$.
\end{remark} 

Note that one can define the $n$-grafting of labelled trees by the $n$-grafting of the trees, where one keeps the labellings on the leaves.

We define the vector space, denoted $Mag^m(V)$, as follows:
$$Mag^m:=\oplus_n  Mag_n^m\otimes V^{\t n},$$
where $Mag^m_n:=\mathbb K [Y_n^m]$.

\begin{proposition}
Let $V$ be a vector space.
The space $Mag^m(V)$ endowed with  the $n$-grafting of labelled trees, for all $2\leq n\leq m$, is a $m$-magmatic algebra. Moreover it is the free $m$-magmatic algebra over $V$.
\end{proposition}$\Box$ 

\subsection{ $m$-magmatic coalgebra and cofree $m$-magmatic coalgebra}

\begin{definition}
An \emph{$m$-magmatic coalgebra} $C$ is a vector space endowed
 with one $n$-ary co-unitary co-operation $\Delta_n:C\to C^{\otimes n}$ for all $2\leq n\leq m$ 
such that:

every $\Delta_n$ admits the same co-unit $c$ and that the following diagram is commutative:
$$
\coproduit
$$
\end{definition}

\subsubsection{The cofree $m$-ary magmatic coalgebra}\label{mcopr}

We endow the vector space of  $m$-planar trees,defined above, with the following $n$-ary co-operations, for $2\leq n\leq m$:

for any $m$-ary planar rooted tree $t$ we define:
$$
\Delta_n(t)=\sum t_1\otimes\cdots \otimes t_n
$$
where the sum is extended to all the manner to write $t=\vee_n(t_1,\cdot,t_n)$.

Explicitely for $t=\vee_n(t_1,\cdots,t_n)$, where all $t_i\neq\emptyset$, we have: 
\begin{eqnarray*}
&&\Delta_n(t):=
\left(\deconcat \right)+\sum_{i=0}^{n-1}\emptyset^{\otimes i}\otimes t\otimes \emptyset^{\otimes n-i-1}\ ,\\
&&\Delta_m(t):=\left\{ \begin{array}{ll}
\sum_{i=0}^{m-1}\emptyset^{\otimes i}\otimes t\otimes \emptyset^{\otimes m-i-1}\ ,
& \textrm{if $m<n$}\\
\\
\sum_{i=0}^{m-1}\emptyset^{\otimes i}\otimes t\otimes \emptyset^{\otimes m-i-1}+\\
+\sum_{i_1+\cdots+i_{n+1}=m-n}\emptyset^{\otimes i_1}\otimes t_1 \otimes \emptyset^{\otimes i_2}\otimes \cdots\otimes t_n\otimes \emptyset^{\otimes i_n+1}
\ ,
& \textrm{if $m>n$}\end{array} \right.\\
&&\Delta_n(|):=\sum_{i=0}^{n-1} \emptyset^{\otimes i}\otimes |\otimes \emptyset^{n-i-1} \ ,\\
&&\Delta_n(\emptyset):=\emptyset^{\otimes n} \ .
\end{eqnarray*}
As in the preceding section one can define the $n$-ungrafting of labelled trees by
the $n$-ungrafting of planar trees and keeping the labelling on the leaves.
 
Remark that $\emptyset$ plays here the role of the unit, it can then be denoted by $1:=\emptyset$.

\begin{proposition}
Let $V$ be a vector space.
The space $Mag^m(V)$ endowed with the $n$-ungrafting co-operations, $2\leq n\leq m$, on labelled trees
 is a connected infinite magmatic coalgebra. Moreover it is free among
over $V$  the connected infinite magmatic coalgebras.
\end{proposition}

\begin{proof} It is similar to the proof of cofree infinite magmatic coalgebra.
\end{proof}

\subsection{ $m$-magmatic bialgebra}

\begin{definition}
An \emph{$m$-magmatic bialgebra} $(\mathcal{H},\mu_n,\Delta_n)$, where $2\leq n\leq m$, is a vector space
$\mathcal{H}=\bar{ \mathcal{H}} \oplus \mathbb{K} 1$ such that:\\
1) $\mathcal{H}$ admits an $m$-magmatic algebra structure with $n$-ary operations denoted~$\mu_n$\!\!~,\\
2) $\mathcal{H}$ admits a $m$-magmatic coalgebra structure with $n$-ary co-operations denoted~$\Delta_n$,\\
3) $\mathcal{H}$ satisfies the following ``compatibility relation'':
\begin{equation}\label{rel m}
\begin{array}{l}
\Delta_n\circ\mu_n (x_1\otimes\cdots\otimes x_n)=
x_1\otimes\cdots\otimes x_n+\sum_{i=0}^{n-1}1^{\otimes i}\otimes \underline x \otimes 1^{\otimes n-i-1}\ ,\\
\Delta_m\circ\mu_n (x_1\otimes\cdots\otimes x_n)=\\
\qquad\qquad\qquad\left\{ \begin{array}{ll}
\sum_{i=0}^{m-1}1^{\otimes i}\otimes \underline x \otimes 1^{\otimes m-i-1}\ ,
& \textrm{if $m<n$}\\
\\
\sum_{i=0}^{m-1}1^{\otimes i}\otimes \underline x \otimes 1^{\otimes m-i-1}+\\
+\sum_{i_1+\cdots+i_n+1=m-n}1^{\otimes i_1}\otimes x_1 \otimes 1^{\otimes i_2}\otimes \cdots\otimes x_n\otimes 1^{\otimes i_{n+1}}
& \textrm{if $m>n$}\end{array} \right.
\end{array}
\end{equation} 
$\forall  \underline x :=\mu_n \circ x_1\otimes\cdots\otimes x_n \ and \  x_1,\cdots,x_n \in \bar{\mathcal{H}}$ and $2\leq n\leq m$.
\end{definition}

\begin{example}Let $V$ be a vector space.
The space
$(Mag^m(V),\vee_n,\Delta_n)$, where $2\leq n\leq m$ and the operations $\vee_n$ (resp. the cooperations $\Delta_n$) as defined in \ref{mpr} and \ref{mcopr} is an
 infinite magmatic connected bialgebra.
\end{example}

\subsection{ The rigidity theorem}

We can now state a rigidity theorem:

\begin{theorem}\label{marymagmatique}
If  $\mathcal{H}$ be a connected $m$-ary magmatic bialgebra over a field $\mathbb{K}$ of any characteristic, then the following are equivalent:
\begin{enumerate}
\item $\mathcal{H}$ is connected,
\item $\mathcal{H} \cong Mag^m(\Prim \mathcal{H})$.
\end{enumerate}
\end{theorem}$\Box$

The proof is very similar to the infinite magmatic case. One as to use the following two definitions and the two lemmas:

\begin{definition}
The \emph{completed $m$-ary magmatic algebra}, denoted by
$Mag^m(\mathbb K)^{\wedge}\ ,$ is defined by
$$Mag^m(\mathbb K)^{\wedge}~\!\!=~\!\!\prod_{n\geq 0}Mag_n \ ,$$ where
 the first generator  $|$ is denoted by $t$. 
This definition allows us to define
 formal power series of trees in $Mag^m(\mathbb K)^{\wedge}$.
\end{definition}

\begin{definition}
Let $2\leq n\leq m$.
The \emph{$n$-convolution} of $n$ $m$-ary magmatic algebra morphisms $f_1,\cdots ,f_n$ is defined by:
\begin{equation*}
 \star_n (f_1\cdots f_n) := \mu_n \circ (f_1 \otimes \cdots \otimes f_n) \circ \Delta_n  \ .
\end{equation*}
Observe that these operations are unitary.
\end{definition}

\begin{lemma}
The following two formal power series are inverse for composition in
  $Mag^m(\mathbb K)^{\wedge}$:
$$g(|):=|-\arbreA-\arbreAt-\arbreAttt-\cdots -\mcorol, \qquad f(|):= \sum T, $$ where the sum is extended to all $m$-ary planar trees $T$. 
\end{lemma}$\Box$

\begin{lemma}
Let $(\mathcal{H}, \mu_n, \Delta_n)$ be a connected $m$-ary  magmatic bialgebra. The linear map
$e:\mathcal{H} \rightarrow \mathcal{H}$ defined as: 
\begin{equation*}
e:=J- \star_2 J^{\otimes 2}-\star_3 J^{\otimes 3}-\cdots-\star_m J^{\otimes m}
\end{equation*}
where $J= Id-uc$, $u$ the unit of the operations, $c$ the co-unit of the co-operations,
has the following properties:
\begin{enumerate}
\item $\Im e = \Prim \mathcal{H} $,
\item forall $x_1,\cdots,x_n \in \bar{ \mathcal{H}}$ one has $ \ e\circ\mu_n (x_1\otimes \cdots \otimes x_n)=0$,
\item the linear map $e$ is an idempotent,
\item for $\mathcal{H}=(Mag^m(V),\mu_n,\Delta_n)$ defined above, $e$ is the identity on $V = Mag_1(V)$
and trivial on the other components.
\end{enumerate}
\end{lemma}$\Box$

\begin{remark}
 Theorems \ref{magmatique}, \ref{marymagmatique}, \ref{theo2} fit into the framework of ``triples of operad'' of Loday (\cite{L3}). In our case the triples is $(Mag^m,Mag^m,Vect)$, $m\in \mathbb N$ or $m=\infty$. In the case $(Mag^m,Mag^n,\mathfrak P)$ the operad of the primitives is still to be unraveled. 
\end{remark}

\section{Binary magmatic bialgebras}

A special case of the $m$-ary magmatic bialgebra is the binary magmatic bialgebras, known also as magmatic bialgebras \cite{G},\cite{H}.

\begin{definition}
A \emph{binary magmatic algebra} $A$ is a   vector space endowed with a binary unitary operation.
\end{definition}

\begin{definition}
A \emph{binary magmatic coalgebra} $A$ is a   vector space endowed with a binary counitary cooperation.
\end{definition}

\subsection{Planar binary trees equipped with a  product and a  coproduct}\label{arbrepl}
A \emph{planar binary tree}  $T$ is a planar connected graph which is assumed to be simple (no loops nor multiple edges), rooted and such that any inner vertices is of valence at most $2$. 
We denote by $Y_n^2$ the set of planar trees with $n$ leaves.
In low dimensions one gets:
\begin{eqnarray*}
&&Y_0^2=\left\{ \emptyset \right\}  ,
\ Y_1^2=\left\{ | \right\} ,
\ Y_2^2=\left\{ \arbreA \right\},
Y_3^2= \left\{ \arbreAB \ \arbreBA  \right\},\\
&&Y_4^2=\left\{ \arbreABC \ \arbreBAC \ \arbreACA \ \arbreCAB \arbreCBA  \right\} , \cdots
\end{eqnarray*}
The \emph{grafting} of two trees is the gluing of the  root of each tree on a new root.
For example, if $t$ and $s$ are two planar binary trees, their grafting is defined as:
$$
\vee_2(t,s):=
\concat \ .
$$

\begin{remark}From our definition of a planar binary tree, any $t\in Y_n^2$ is of the form
$$ t=\vee (t_1,t_2)$$
for uniquely determined trees $t_1,t_2$. Diagrammatically it is represented by:
$$
t= \concattt
$$
\end{remark}  

Thanks to this remark one can define the un-grafting of a tree $t = t_1 \vee t_2$ as:

$$
\Delta(t):= \deconcattwo +t\otimes \emptyset+\emptyset\otimes t\ .
$$
We add: 
\begin{eqnarray*}
\Delta (|)&=&|\otimes\emptyset+\emptyset\otimes | \ , \\
\Delta (\emptyset) &=&\emptyset\otimes \emptyset \ .
\end{eqnarray*}

Let $V$ be a vector space, one can define the grafting (resp. the ungrafting) of labelled trees by the grafting (resp. the ungrafting) of the trees, where one keeps the labellings on the leaves.
Therefore we have defined an operation and a cooperation on the vector space $Mag^2(V):=\oplus_{n\geq 0}Mag_n^2\otimes V^{\otimes n}$, where $Mag_n^2=\mathbb K [Y_n^2]$ .

\subsection{Binary magmatic bialgebras}

\begin{definition}
A \emph{binary magmatic bialgebra} $(\mathcal{H},\cdot,\Delta)$ is a vector space such that 
$\mathcal{H}=\bar{ \mathcal{H}} \oplus \mathbb{K}\cdot 1$ verifying:\\
1) $(\mathcal{H},\cdot)$ is a binary magmatic algebra,\\
2) $(\mathcal{H},\Delta)$ is a binary magmatic coalgebra,\\
3) $\mathcal{H}$ satisfies the following ``compatibility relation'' :
\begin{equation*}
\Delta(x \cdot y)=x \cdot y \otimes 1 + x \otimes y + 1 \otimes x \cdot y, \ \forall  x,y \in \mathcal{H} \ .
\end{equation*}
\end{definition}

\begin{definition}
A binary magmatic coalgebra $\mathcal{H}$ is said to be \emph{connected} if it verifies the following property:
\begin{equation*}
\begin{array}{l}
\mathcal{H}= \bigcup_{r \geq 0}F_{r}\mathcal{H}\qquad
\textrm{where } F_{0}:=\mathbb{K} 1\\
\textrm{and, by induction } F_{r}:=\left\{x \in \mathcal{H} \ \vert \  \bar{ \Delta}(x) \in F_{r-1} \otimes F_{r-1}\right\} \ .
\end{array}
\end{equation*}
where $\bar{\Delta}(x)=\Delta(x)-x\otimes 1-1 \otimes x$ .\\
Note that the connectedness only depends on the unit and the cooperation.
\end{definition}

A fundamental example, in this context, is:
\begin{proposition}
The triple $(Mag^2(V),\cdot,\Delta)$, where $\cdot$ (resp. $\Delta$) were defined in \ref{arbrepl} is a connected magmatic bialgebra.
\end{proposition}

And, we can sate the following theorem:
\begin{theorem}\label{theo2}
Let $\mathcal{H}$ be a connected binary magmatic bialgebra over a field of any caracteristic, then $\mathcal{H}$ is isomorphic to the binary magmatic bialgebra $Mag^2(\Prim \mathcal{H})$.
\end{theorem}

\end{document}